\documentclass{article}
\usepackage{amsmath,amssymb,url,amsthm,cases}
\usepackage[pdftex,hiresbb]{graphicx}
\usepackage{type1cm}
\usepackage{subcaption}
\usepackage{geometry}
\usepackage{here}
\allowdisplaybreaks
\newtheorem{thm}{Theorem}[section]

\newtheorem{lem}[thm]{Lemma}
\newtheorem{prop}[thm]{Proposition}

\newtheorem{rem}[thm]{Remark}

\newtheorem{prob}[thm]{Problem}
\def\l     {\left}
\def\r     {\right}
\def\<     {\langle}
\def\>     {\rangle}
\def\fin   {\hfill{$\Box$}\vspace{5mm}}
\def\bbE   {{\mathbb E}}
\def\bbF   {{\mathbb F}}
\def\bbG   {{\mathbb G}}
\def\bbN   {{\mathbb N}}
\def\bbP   {{\mathbb P}}
\def\bbR   {{\mathbb R}}

\def\calA  {{\mathcal A}}
\def\calF  {{\mathcal F}}
\def\calG  {{\mathcal G}}
\def\calP  {{\mathcal P}}
\def\calQ  {{\mathcal Q}}
\def\calT  {{\mathcal T}}

\def\ve    {\varepsilon}
\def\olA   {\overline{A}}
\def\olV   {\overline{V}}
\def\olS   {\overline{S}}
\def\olZ   {\overline{Z}}

\def\tK   {\widetilde{K}}
\def\tT   {\widetilde{T}}
\def\tpi   {\widetilde{\pi}}
\def\hT   {\widehat{T}}
\def\olx   {\widehat{x}}
\def\oli   {\overline{i}}
\def\olx   {\overline{x}}
\begin{document}
\title{Constrained optimal stopping under a regime-switching model}
%\author{Takuji Arai\footnote{Department of Economics, Keio University, 2-15-45 Mita, Minato-ku, Tokyo, 108-8345, Japan. \\
% (arai@econ.keio.ac.jp)} \\
%Masahiko Takenaka\footnote{Department of Economics, Keio University, 2-15-45 Mita, Minato-ku, Tokyo, 108-8345, Japan.}}
\author{Takuji Arai and Masahiko Takenaka \\ Department of Economics, Keio University}
\maketitle

%%%%%%%%%%%%%%%%%%%%%%%%%%%%%%%%%%%%%%%%%%%%%%%%%%%%%%%%%%%%%%%%%%%%%%%%%%%%%%%
\begin{abstract}
We investigate an optimal stopping problem for the expected value of a discounted payoff on a regime-switching geometric Brownian motion
under two constraints on the possible stopping times: only at exogenous random times and only during a specific regime.
The main objectives are to show that an optimal stopping time exists as a threshold type under some boundary conditions and
to derive expressions of the value functions and the optimal threshold.
To this end, we solve the corresponding variational inequality and show that its solution coincides with the value functions.
Some numerical results are also introduced. Furthermore, we investigate some asymptotic behaviors. 
\end{abstract}

%%%%%%%%%%%%%%%%%%%%%%%%%%%%%%%%%%%%%%%%%%%%%%%%%%%%%%%%%%%%%%%%%%%%%%%%%%%%%%%
%
% Section 1
%
%%%%%%%%%%%%%%%%%%%%%%%%%%%%%%%%%%%%%%%%%%%%%%%%%%%%%%%%%%%%%%%%%%%%%%%%%%%%%%%
\setcounter{equation}{0}
\section{Introduction}
In the real options literature, the following type of optimal stopping problems appears frequently:
\begin{equation}\label{primal-prob}
\sup_{\tau\in\calT}\bbE\l[e^{-r\tau}\pi(X_\tau)\big|X_0=x\r],
\end{equation}
where $r>0$ is the exogenous discount rate, $X=\{X_t\}_{t\geq0}$ is a stochastic process, which we call the cash flow process,
$\calT$ is the set of all stopping times that investors can choose, and $\pi$ is an $\bbR$-valued function, which we call the payoff function.
We can regard (\ref{primal-prob}) as a function on $x$, which we call the value function.
Problem (\ref{primal-prob}) concerns the optimal investment timing for an investment whose payoff is given by the random variable $\pi(X_t)$ when executed at time $t$.
The most typical example of $\pi$ is
\begin{equation}\label{primal-pi}
\pi(x)=\bbE\l[\int_t^\infty e^{-r(s-t)}X_sds-I \ \Big| \ X_t=x\r],
\end{equation}
which expresses the value of an investment that starts at time $t$ with an initial cost $I>0$ and that brings to the investor perpetually an instantaneous return $X_s$ at each time $s>t$.
Remark that the right-hand side of (\ref{primal-pi}) becomes a function on $x$ when the process $X$ has the strong Markov property such as a geometric Brownian motion.
The main concern of (\ref{primal-prob}) is to show that an optimal stopping time $\tau^*\in\calT$ exists and can be expressed as
\[
\tau^*=\inf\{t>0|X_t\geq x^*\}
\]
for some $x^*\in\bbR$. This type of optimal stopping is called threshold type, and $x^*$ is called its optimal threshold.
It is significant to examine whether an optimal stopping is of threshold type.
If so, the optimal strategy becomes apparent, and the optimal stopping time can be explicitly described.
McDonald and Siegel \cite{MS86} has undertaken this framework of optimal stopping problems. See also Chapter 5 of Dixit and Pindyck \cite{DP12}.
Here we focus on discussing (\ref{primal-prob}) when $X_t$ is a regime-switching geometric Brownian motion under two constraints on $\calT$.

Regime-switching models, widely studied in mathematical finance (\cite{B98}, \cite{BE02}, \cite{BE02b}, \cite{ECS05}, \cite{G01}, \cite{GZ04} and so forth),
are models in which the regime, representing, e.g., the economy's general state, changes randomly.
In this paper, we consider a regime-switching model with two regimes $\{0,1\}$. Let $\theta=\{\theta_t\}_{t\geq0}$ be a stochastic process expressing the regime at time $t$.
In particular, $\theta$ is a $\{0,1\}$-valued continuous-time Markov chain.
Then the cash-flow process $X$ is given by the solution to the following stochastic differential equation (SDE):
\begin{equation}\label{X-SDE}
dX_t=X_t(\mu_{\theta_t}dt+\sigma_{\theta_t}dW_t), \ \ \ X_0>0,
\end{equation}
where $\mu_i\in\bbR$ and $\sigma_i>0$ for $i=0,1$, and $W=\{W_t\}_{t\geq0}$ is a one-dimensional standard Brownian motion independent of $\theta$.
Considering a regime-switching model, we need to define a value function for each initial regime, that is, for each $i=0,1$, we define the value function $v_i$ as
\begin{equation}\label{def-vi-Sect1}
v_i(x):=\sup_{\tau\in\calT}\bbE\l[e^{-r\tau}\pi(X_\tau)\big|\theta_0=i,X_0=x\r].
\end{equation}

Furthermore, we impose two constraints on $\calT$.
Liquidity risk and other considerations mean that investment is not always possible. Therefore, it is significant to analyze models with constraints on investment opportunities and timing.
Hence, we impose two constraints simultaneously in this paper. One is the random arrival of investment opportunities.
More precisely, we restrict stopping only at exogenous random times given by the jump times of a Poisson process independent of $W$ and $\theta$.
Another is the regime constraint. We restrict that stopping is feasible only during regime $1$.

Now, we introduce some related works.
Bensoussan et al. \cite{BYY12} discussed the problem (\ref{primal-prob}) for the same cash flow process $X$ as defined in (\ref{X-SDE}) without restriction on stopping.
They treated the case where $\pi$ is given as (\ref{primal-pi}) and showed that an optimal stopping time exists as a threshold type by an argument based on PDE techniques.
Nishihara \cite{N20} discussed the same problem for a two-state regime-switching model with $\pi(x)=x-I$ under the regime constraint,
but the cash flow process is still a geometric Brownian motion. Note that \cite{N20} assumed that an optimal stopping exists as a threshold type.
In addition, Egami and Kevkhishvili \cite{EK20} also studied the same problem for the case where $X$ is a regime-switching diffusion process but without restriction on stopping.
On the other hand, the restriction of stopping to exogenous random times has been undertaken by Dupuis and Wang \cite{DW02}.
They considered the case where the cash flow process is a geometric Brownian motion and the payoff function is of American call option type, i.e., $\pi(x)=(x-K)^+$,
and did not deal with regime-switching models.
In \cite{DW02}, they first derived a variational inequality (VI) through a heuristic discussion.
Solving it, they showed by a probabilistic argument that the solution to the VI coincides with the value function.
There are other many works dealing with this issue such as \cite{H21}, \cite{HZ19}, \cite{LRS20}, \cite{L12}, \cite{MR16} and so forth.

To our best knowledge, this paper is the first study that deals with the constrained optimal stopping problem on a regime-switching geometric Brownian motion.
It is also new to simultaneously impose the random arrival of investment opportunities and the regime constraint.
Remark that the discussion in this paper is based on the approach in \cite{DW02}.

This paper is organized as follows: Some mathematical preparations and the formulation of our optimal stopping problem will be given in Section 2.
Section 3 introduces the corresponding VI and solves its modified version in which two boundary conditions are replaced.
We shall derive explicit expressions of the solution to the modified VI, which involves solutions to quartic equations, but it can be numerically computable easily.
In Section 4, assuming that the two boundary conditions replaced in Section 3 are satisfied,
we prove that the solution to the VI coincides with the value functions and the optimal threshold for our optimal stopping problem.
In addition, we introduce some numerical results. Section 5 is devoted to illustrating some results on asymptotic behaviors, and Section 6 concludes this paper.

%%%%%%%%%%%%%%%%%%%%%%%%%%%%%%%%%%%%%%%%%%%%%%%%%%%%%%%%%%%%%%%%%%%%%%%%%%%%%%%
%
% Section 2
%
%%%%%%%%%%%%%%%%%%%%%%%%%%%%%%%%%%%%%%%%%%%%%%%%%%%%%%%%%%%%%%%%%%%%%%%%%%%%%%%
\setcounter{equation}{0}
\section{Preliminaries and problem formulation}
We consider a regime-switching model with state space $\{0,1\}$ and suppose that the regime process $\theta$ is a $\{0,1\}$-valued continuous-time Markov chain with generator
\[\begin{pmatrix}
-\lambda_0 & \lambda_0 \\
\lambda_1  & -\lambda_1
\end{pmatrix},\]
where $\lambda_0,\lambda_1>0$. Now, we make the convention $\theta_\infty\equiv1$. Note that the length of regime $i$ follows the exponential distribution with parameter $\lambda_i$.
We take the process $X$ defined in (\ref{X-SDE}) as the cash flow process, and assume throughout this paper that
\begin{equation}\label{r-mu}
r>\mu_0\vee\mu_1.
\end{equation}
Let $J=\{J_t\}_{t\geq0}$ be a Poisson process with intensity $\eta>0$ independent of $W$ and $\theta$,
and denote by $T_k$ its $k$th jump time for $k\in\bbN$ with the conventions $T_0\equiv0$ and $T_\infty\equiv\infty$, where $\bbN:=\{1,2,\dots\}$.
Note that the process $J$ generates exogenous random times when an investment opportunity arrives.
In other words, for $k\in\bbN$, $T_k$ represents the $k$th investment opportunity time.
Suppose that $\theta$, $W$, and $J$ are defined on a complete probability space $(\Omega,\calF,\bbP)$.
In addition, we denote by $\bbF=\{\calF_t\}_{t\geq0}$ the filtration generated by $\theta$, $W$ and $J$. Assume that $\bbF$ satisfies the usual condition.
Furthermore, we restrict stopping to only when the regime is $1$. Thus, the set of all possible stopping times is described by
\[
\calT:=\{\tau\in\calT_0 \ | \ \mbox{ for each }\omega\in\Omega, \theta_{\tau(\omega)}(\omega)=1\mbox{ and }\tau(\omega)=T_j(\omega)\mbox{ for some }j\in\bbN_\infty\},
\]
where $\calT_0$ is the set of all $[0,\infty]$-valued stopping times and $\bbN_\infty:=\bbN\cup\{\infty\}$.
Next we formulate the payoff function $\pi$ as follows:
\begin{equation}\label{eq-pi}
\pi(x)=\alpha(x-K)^+-I
\end{equation}
for some $\alpha>0$, $K\geq0$ and $I\geq0$, but we exclude the case where $K=I=0$ since the optimal threshold $x^*$ is obviously $0$ in this case. 
This formulation includes $\pi(x)=(x-K)^+$ treated in \cite{DW02}, and $\pi(x)=x-I$ in \cite{N20}.
Moreover, (\ref{eq-pi}) covers the payoff function introduced in (\ref{primal-pi}).
In fact, \cite{BYY12} showed that
\[
\bbE\l[\int_0^\infty e^{-rt}X_tdt\Big|\theta_0=i,X_0=x\r]=\frac{(r-\mu_{1-i}+\lambda_i+\lambda_{1-i})x}{(r-\mu_{1-i})(r-\mu_i)+\lambda_i(r-\mu_{1-i})+\lambda_{1-i}(r-\mu_i)}.
\]

In the setting described above, we define the value functions $v_i, i=0,1$ as follows:
\begin{equation}\label{def-vi}
\l\{\begin{array}{l}
v_1(x):=\displaystyle{\sup_{\tau\in\calT}}\bbE^{1,x}\l[e^{-r\tau}\pi(X_\tau)\r], \vspace{1.5mm} \\
v_0(x):=\bbE^{0,x}\l[e^{-r\xi_0}v_1(X_{\xi_0})\r]
\end{array}\r.
\end{equation}
for $x>0$, where $\xi_0:=\inf\{t>0|\theta_t=1\}$ and $\bbE^{i,x}$ means the expectation with the initial condition $\theta_0=i$ and $X_0=x$.
In fact, we should define $v_0$ as $v_0(x):=\sup_{\tau\in\calT}\bbE^{0,x}\l[e^{-r\tau}\pi(X_\tau)\r]$ in terms of (\ref{def-vi-Sect1}),
but the above definition (\ref{def-vi}) is justified by the following:
\begin{align*}
\sup_{\tau\in\calT}\bbE^{0,x}\l[e^{-r\tau}\pi(X_\tau)\r]
=\bbE^{0,x}\l[e^{-r\xi_0}\sup_{\tau^\prime\in\calT^\prime}\bbE\l[e^{-r\tau^\prime}\pi(X^\prime_\tau)|\theta^\prime_0=1,X^\prime_0=X_{\xi_0}\r]\r]
=\bbE^{0,x}\l[e^{-r\xi_0}v_1(X_{\xi_0})\r],
\end{align*}
where $\theta^\prime$ and $X^\prime$ are independent copies of $\theta$ and $X$, respectively,
and $\calT^\prime$ is the set of all possible stopping times defined based on $\theta^\prime$ and $X^\prime$.
We discuss the optimal stopping problem (\ref{def-vi}) in the following sections.

%%%%%%%%%%%%%%%%%%%%%%%%%%%%%%%%%%%%%%%%%%%%%%%%%%%%%%%%%%%%%%%%%%%%%%%%%%%%%%%
%
% Section 3
%
%%%%%%%%%%%%%%%%%%%%%%%%%%%%%%%%%%%%%%%%%%%%%%%%%%%%%%%%%%%%%%%%%%%%%%%%%%%%%%%
\setcounter{equation}{0}
\section{Variational inequality}
We discuss the variational inequality (VI) corresponding to the value functions $v_i, i=0,1$.
From the same sort of argument as Section 3 in \cite{DW02}, the VI is given as follows:

%%%%%%%%%%%%%%%%%%%%%%%%%%%%%%%%%%%%%%%%%%%%%%%%%%%%%%%%%%%%%%%%%%%%%%%%%%%%%%%
\begin{prob}\label{prob-VI}
Find two nonnegative $C^2$-functions $V_0, V_1:\bbR_+\to\bbR_+$ and a constant $x^*\geq\tK$ satisfying
\begin{numcases}{}
V_i(0+)= 0, \ \ \ i=0,1, \label{VI-0} \\
-rV_0(x)+\calA_0 V_0(x)+\lambda_0(V_1(x)-V_0(x))= 0, \ \ \ x>0, \label{VI-ODE0} \\
-rV_1(x)+\calA_1 V_1(x)+\lambda_1(V_0(x)-V_1(x))= 0, \ \ \ 0<x<x^*, \label{VI-ODE1-1} \\
-rV_1(x)+\calA_1 V_1(x)+\lambda_1(V_0(x)-V_1(x))+\eta(\pi(x)-V_1(x))= 0, \ \ \ x>x^*, \label{VI-ODE1-2} \\
V_1(x^*)=\pi(x^*), \label{VI-*} \\
V_1(x)>\pi(x), \ \ \ 0<x<x^*, \label{VI-*L} \\
V_1(x)<\pi(x), \ \ \ x>x^*, \label{VI-*U} 
\end{numcases}
where $\bbR_+:=[0,\infty)$, $\tK:=K+\displaystyle{\frac{I}{\alpha}}$, and $\calA_i, i=0,1$ are the infinitesimal generators of $X$ under regime $i$ defined as
\[
(\calA_if)(x):=\mu_ixf^{\prime}(x)+\frac{1}{2}\sigma_i^2x^2f^{\prime\prime}(x), \ \ \ x>0,
\]
for $C^2$-function $f$.
\end{prob}

This section aims to solve the following modified version of Problem \ref{prob-VI},
in which we replace the boundary conditions (\ref{VI-*L}) and (\ref{VI-*U}) with (\ref{preVI-*U}) below:

%%%%%%%%%%%%%%%%%%%%%%%%%%%%%%%%%%%%%%%%%%%%%%%%%%%%%%%%%%%%%%%%%%%%%%%%%%%%%%%
\begin{prob}\label{prob-preVI}
Find two $C^2$-functions $V_0, V_1:\bbR_+\to\bbR_+$ and a constant $x^*\geq\tK$ satisfying (\ref{VI-0}) -- (\ref{VI-*}) and
\begin{equation}\label{preVI-*U}
0<\lim_{x\to\infty}\frac{V_1(x)}{\pi(x)}<1. 
\end{equation}
\end{prob}

\noindent
To solve Problem \ref{prob-preVI}, we need some preparations.
For $i=0,1$ and $k=L,U$, $G^k_i$ is the quadratic function on $\beta\in\bbR$ defined as
\[
G^k_i(\beta):=\frac{1}{2}\sigma^2_i\beta(\beta-1)+\mu_i\beta-(\lambda_i+r+\eta{\bf 1}_{\{i=1,k=U\}}).
\]
The equation $G^k_i(\beta)=0$ has one positive and one negative solution, denoted by $\zeta^{k,+}_i$ and $\zeta^{k,-}_i$, respectively.
For each $k=U,L$, we denote
\[
F^k(\beta):=G^k_0(\beta)G^k_1(\beta)-\lambda_0\lambda_1,
\]
and consider the quartic equation $F^k(\beta)=0$. Since $F^k(0)>0$, $F^k(\zeta^{k,\pm}_i)<0$, and $F^k(\beta)\to\infty$ as $\beta$ tends to $\pm\infty$,
the equation $F^k(\beta)=0$ has four different solutions, two of which are positive, and two of which are negative.
Now, for the equation $F^L(\beta)=0$, we denote the larger positive solution by $\beta^L_A$ and another positive solution by $\beta^L_B$.
Note that $F^L(1)$ is positive, and $\zeta^{L,+}_i>1$ holds since $G^L_i(1)<0$. Thus, $1<\beta^L_B<\zeta^{L,+}_i<\beta^L_A$ holds for $i=0,1$.
A similar argument can be found in Remark 2.1 of Guo \cite{G01}.
Furthermore, the same holds for the quartic equation $F^U(\beta)=0$.
Let $\beta^U_A$ and $\beta^U_B$ be the larger and other negative solutions to $F^U(\beta)=0$, respectively, that is, $\beta^U_B<\zeta^{U,-}_i<\beta^U_A<0$ holds for $i=0,1$.
In addition, we define the following constants:
\begin{equation}\label{eq-singular}
\l\{\begin{array}{l}
a_0:= \displaystyle{\frac{\alpha\eta\lambda_0}{(r-\mu_0+\lambda_0)(r-\mu_1+\lambda_1+\eta)-\lambda_0\lambda_1}}, \ \ \
a_1:= \displaystyle{\frac{\alpha\eta(r-\mu_0+\lambda_0)}{(r-\mu_0+\lambda_0)(r-\mu_1+\lambda_1+\eta)-\lambda_0\lambda_1}}, \vspace{1.5mm} \\
b_0:= \displaystyle{\frac{\alpha\tK\eta\lambda_0}{\lambda_0\lambda_1-(r+\lambda_0)(r+\lambda_1+\eta)}}, \ \ \
b_1:= \displaystyle{\frac{\alpha\tK\eta(r+\lambda_0)}{\lambda_0\lambda_1-(r+\lambda_0)(r+\lambda_1+\eta)}},
\end{array}\r.
\end{equation}
and
\begin{equation}\label{eq-PQ}
\l\{\begin{array}{l}
P^L_A:= \displaystyle{\frac{\alpha(\beta^L_B-\beta^U_A)(-\beta^L_B+\beta^U_B)+a_1(\beta^U_A-1)(\beta^U_B-1)}{(\beta^L_A-\beta^L_B)(\beta^L_A+\beta^L_B-\beta^U_A-\beta^U_B)}}, \vspace{1.5mm} \\
Q^L_A:= \displaystyle{\frac{-\alpha\tK(\beta^L_B-\beta^U_A)(-\beta^L_B+\beta^U_B)+b_1\beta^U_A\beta^U_B}{(\beta^L_A-\beta^L_B)(\beta^L_A+\beta^L_B-\beta^U_A-\beta^U_B)}}, \vspace{1.5mm} \\
P^L_B:= \displaystyle{\frac{\alpha(\beta^L_A-\beta^U_A)(\beta^L_A-\beta^U_B)-a_1(\beta^U_A-1)(\beta^U_B-1)}{(\beta^L_A-\beta^L_B)(\beta^L_A+\beta^L_B-\beta^U_A-\beta^U_B)}}, \vspace{1.5mm} \\
Q^L_B:= \displaystyle{\frac{-\alpha\tK(\beta^L_A-\beta^U_A)(\beta^L_A-\beta^U_B)-b_1\beta^U_A\beta^U_B}{(\beta^L_A-\beta^L_B)(\beta^L_A+\beta^L_B-\beta^U_A-\beta^U_B)}}, \vspace{1.5mm} \\
P^U_A:= \displaystyle{\frac{\alpha(\beta^L_A-\beta^U_B)(\beta^L_B-\beta^U_B)+a_1(\beta^U_B-1)(\beta^L_A+\beta^L_B-\beta^U_B-1)}{(\beta^U_A-\beta^U_B)(\beta^L_A+\beta^L_B-\beta^U_A-\beta^U_B)}}, \vspace{1.5mm} \\
Q^U_A:= \displaystyle{\frac{-\alpha\tK(\beta^L_A-\beta^U_B)(\beta^L_B-\beta^U_B)+b_1\beta^U_B(\beta^L_A+\beta^L_B-\beta^U_B)}{(\beta^U_A-\beta^U_B)(\beta^L_A+\beta^L_B-\beta^U_A-\beta^U_B)}}, \vspace{1.5mm} \\
P^U_B:= \displaystyle{\frac{\alpha(\beta^L_A-\beta^U_A)(-\beta^L_B+\beta^U_A)-a_1(\beta^U_A-1)(\beta^L_A+\beta^L_B-\beta^U_A-1)}{(\beta^U_A-\beta^U_B)(\beta^L_A+\beta^L_B-\beta^U_A-\beta^U_B)}}, \vspace{1.5mm} \\
Q^U_B:= \displaystyle{\frac{-\alpha\tK(\beta^L_A-\beta^U_A)(-\beta^L_B+\beta^U_A)-b_1\beta^U_A(\beta^L_A+\beta^L_B-\beta^U_A)}{(\beta^U_A-\beta^U_B)(\beta^L_A+\beta^L_B-\beta^U_A-\beta^U_B)}}.
\end{array}\r.
\end{equation}

With the above preparations, we solve Problem \ref{prob-preVI} as follows:

%\[
%\l\{\begin{array}{l}
%A^L_1(x^*)^{\beta^L_A}= \displaystyle{\frac{(\beta^L_B-\beta^U_A)(-\beta^L_B+\beta^U_B)\pi(x^*)+(\beta^U_A-1)(\beta^U_B-1)a_1x_*+\beta^U_A\beta^U_Bb_1}{(\beta^L_A-\beta^L_B)(\beta^L_A+\beta^L_B-\beta^U_A-\beta^U_B)}}, \vspace{1.5mm} \\
%B^L_1(x^*)^{\beta^L_B}= \displaystyle{\frac{(\beta^L_A-\beta^U_A)(\beta^L_A-\beta^U_B)\pi(x^*)-(\beta^U_A-1)(\beta^U_B-1)a_1x^*-\beta^U_A\beta^U_Bb_1}{(\beta^L_A-\beta^L_B)(\beta^L_A+\beta^L_B-\beta^U_A-\beta^U_B)}}, \vspace{1.5mm} \\
%A^U_1(x^*)^{\beta^U_A}= \displaystyle{\frac{(\beta^L_A-\beta^U_B)(\beta^L_B-\beta^U_B)\pi(x^*)+(\beta^U_B-1)(\beta^L_A+\beta^L_B-\beta^U_B-1)a_1x^*+\beta^U_B(\beta^L_A+\beta^L_B-\beta^U_B)b_1}{(\beta^U_A-\beta^U_B)(\beta^L_A+\beta^L_B-\beta^U_A-\beta^U_B)}}, \vspace{1.5mm} \\
%B^U_1(x^*)^{\beta^U_B}= \displaystyle{\frac{(\beta^L_A-\beta^U_A)(-\beta^L_B+\beta^U_A)\pi(x^*)-(\beta^U_A-1)(\beta^L_A+\beta^L_B-\beta^U_A-1)a_1x^*-\beta^U_A(\beta^L_A+\beta^L_B-\beta^U_A)b_1}{(\beta^U_A-\beta^U_B)(\beta^L_A+\beta^L_B-\beta^U_A-\beta^U_B)}}.
%\end{array}\r.
%\]

%%%%%%%%%%%%%%%%%%%%%%%%%%%%%%%%%%%%%%%%%%%%%%%%%%%%%%%%%%%%%%%%%%%%%%%%%%%%%%%
\begin{prop}\label{prop-preVI}
Problem \ref{prob-preVI} has the following unique solution $(V_0,V_1,x^*)$: For $i=0,1$,
\begin{numcases}{}
V_i(x) = A^L_ix^{\beta^L_A}+B^L_ix^{\beta^L_B}, \ \ \ 0<x<x^*, \label{eq-ViL} \\
V_i(x) = A^U_ix^{\beta^U_A}+B^U_ix^{\beta^U_B}+a_ix+b_i, \ \ \ x>x^*, \label{eq-ViU}
\end{numcases}
and
\begin{equation}\label{eq-x*}
x^*=-\frac{\frac{(1-\beta^L_A)Q^L_A}{G^L_0(\beta^L_A)}+\frac{(1-\beta^L_B)Q^L_B}{G^L_0(\beta^L_B)}+\frac{(\beta^U_A-1)Q^U_A}{G^U_0(\beta^U_A)}
      +\frac{(\beta^U_B-1)Q^U_B}{G^U_0(\beta^U_B)}+\frac{b_0}{\lambda_0}}
      {\frac{(1-\beta^L_A)P^L_A}{G^L_0(\beta^L_A)}+\frac{(1-\beta^L_B)P^L_B}{G^L_0(\beta^L_B)}+\frac{(\beta^U_A-1)P^U_A}{G^U_0(\beta^U_A)}
      +\frac{(\beta^U_B-1)P^U_B}{G^U_0(\beta^U_B)}}
\end{equation}
where
\begin{equation}\label{eq-AB1}
\l\{\begin{array}{l}
A^k_1=(x^*)^{-\beta^k_A}(P^k_Ax^*+Q^k_A), \ \ \ B^k_1=(x^*)^{-\beta^k_B}(P^k_Bx^*+Q^k_B), \vspace{1.5mm} \\
A^k_0=\displaystyle{\frac{-\lambda_0}{G^k_0(\beta^k_A)}A^k_1}, \ \ \ B^k_0=\displaystyle{\frac{-\lambda_0}{G^k_0(\beta^k_B)}B^k_1}
\end{array}\r.
\end{equation}
for $k=L,U$.
\end{prop}

\proof For the time being, we use $\tpi(x):=\alpha x-\alpha\tK$ instead of $\pi$, that is, we rewrite (\ref{VI-ODE1-2}) and (\ref{VI-*}) as follows:
\begin{numcases}{}
-rV_1(x)+\calA_1 V_1(x)+\lambda_1(V_0(x)-V_1(x))+\eta(\tpi(x)-V_1(x))= 0, \ \ \ x>x^*, \label{preVI-ODE1-2} \\
V_1(x^*)=\tpi(x^*). \label{preVI-*}
\end{numcases}
\noindent{\bf Step 1:} \ For $0<x<x^*$, a general solution to (\ref{VI-ODE0}) and (\ref{VI-ODE1-1}) is expressed as (\ref{eq-ViL})
with some $A^L_i,B^L_i\in\bbR$ and some $\beta^L_A,\beta^L_B>0$. Remark that the non-negativity of $\beta^L_A$ and $\beta^L_B$ is derived from the condition (\ref{VI-0}).
Without loss of generality, we may assume that $\beta^L_A>\beta^L_B$.
Substituting (\ref{eq-ViL}) for (\ref{VI-ODE0}) and (\ref{VI-ODE1-1}), we obtain that
\[
\l(A^L_iG^L_i(\beta^L_A)+\lambda_iA^L_{1-i}\r)x^{\beta^L_A}+\l(B^L_iG^L_i(\beta^L_B)+\lambda_iB^L_{1-i}\r)x^{\beta^L_B}=0, \ \ \ i=0,1,
\]
for any $x\in(0,x^*)$, which is equivalent to that $A^L_iG^L_i(\beta^L_A)+\lambda_iA^L_{1-i}=0$ and $B^L_iG^L_i(\beta^L_B)+\lambda_iB^L_{1-i}=0$ for $i=0,1$.
Thus, $\beta^L_A$ satisfies $A^L_0G^L_0(\beta^L_A)A^L_1G^L_1(\beta^L_A)=(-\lambda_0A^L_1)(-\lambda_1A^L_0)$, that is, $G^L_0(\beta^L_A)G^L_1(\beta^L_A)-\lambda_0\lambda_1=0$.
In addition, the same is true for $\beta^L_B$.
Thus, as defined above, $\beta^L_A$ and $\beta^L_B$ are the larger and smaller positive solutions to the equation $F^L(\beta)=0$.
Moreover, $A^L_i$ and $B^L_i$ satisfy the following:
\begin{equation}\label{eq-ABL}
A^L_0=-\frac{\lambda_0}{G^L_0(\beta^L_A)}A^L_1, \mbox{ and } \ B^L_0=-\frac{\lambda_0}{G^L_0(\beta^L_B)}B^L_1.
\end{equation}

%%%%%%%%%%%%%%%%%%%%%%%%%%%%%%%%%%%%%%%%%%%%%%%%%%%%%%%%%%%%%%%%%%%%%%%%%%%%%%%
\noindent{\bf Step 2:} \ Next, we discuss the case where $x>x^*$.
Firstly, we need to find a special solution to (\ref{VI-ODE0}) and (\ref{preVI-ODE1-2}), since (\ref{preVI-ODE1-2}) is inhomogeneous.
Note that $\tpi$ is of linear growth. For each $i=0,1$, we can then write a special solution as $a_ix+b_i$.
Substituting $a_ix+b_i$ for (\ref{VI-ODE0}) and (\ref{preVI-ODE1-2}), we have that
\begin{equation}\label{eq-singular-1}
\l\{\begin{array}{l}
\l(-ra_0+\mu_0a_0+\lambda_0(a_1-a_0)\r)x+\l(-rb_0+\lambda_0(b_1-b_0)\r)=0, \\
\l(-ra_1+\mu_1a_1+\lambda_1(a_0-a_1)+\eta(\alpha-a_1)\r)x+\l(-rb_1+\lambda_1(b_0-b_1)+\eta(-\alpha\tK-b_1)\r)=0
\end{array}\r.
\end{equation}
for any $x>x^*$, in other words, all coefficients in (\ref{eq-singular-1}) are 0, from which $a_i$ and $b_i$ satisfy (\ref{eq-singular}).

Now, we derive $V_i(x)$ for $x>x^*$ in the same way as the previous step.
For each $i=0,1$, we can write a general solution to (\ref{VI-ODE0}) and (\ref{preVI-ODE1-2}) as
\[
V_i(x)=A^U_ix^{\beta^U_A}+B^U_ix^{\beta^U_B}+a_ix+b_i, \ \ \ x>x^*
\]
with some $A^U_i,B^U_i\in\bbR$ and $\beta^U_A,\beta^U_B\in\bbR$.
By (\ref{VI-ODE0}), (\ref{preVI-ODE1-2}) and (\ref{eq-singular-1}), it follows that
\[
A^U_iG^U_i(\beta^U_A)+\lambda_iA^U_{1-i}=0, \mbox{ and } \ B^U_iG^U_i(\beta^U_B)+\lambda_iB^U_{1-i}=0
\]
for $i=0,1$. Thus, by the same way as Step 1, $\beta^U_A$ and $\beta^U_B$ are solutions to the quartic equation $F^U(\beta)=0$.
On the other hand, if either at least $\beta^U_A$ or $\beta^U_B$ is positive, then (\ref{preVI-*U}) is violated since any positive solution is greater than 1.
Thus, $\beta^U_A$ and $\beta^U_B$ are the negative solutions, and we may take them so that $\beta^U_B<\beta^U_A<0$ without loss of generality.
Moreover, we have
\begin{equation}\label{eq-ABU}
A^U_0=-\frac{\lambda_0}{G^U_0(\beta^U_A)}A^U_1, \mbox{ and } \ B^U_0=-\frac{\lambda_0}{G^U_0(\beta^U_B)}B^U_1.
\end{equation}

%%%%%%%%%%%%%%%%%%%%%%%%%%%%%%%%%%%%%%%%%%%%%%%%%%%%%%%%%%%%%%%%%%%%%%%%%%%%%%%
\noindent{\bf Step 3:} \ By the $C^2$-property of $V_1$ and the boundary condition (\ref{preVI-*}), it follows that
\begin{numcases}{}
A^L_1(x^*)^{\beta^L_A}+B^L_1(x^*)^{\beta^L_B}=A^U_1(x^*)^{\beta^U_A}+B^U_1(x^*)^{\beta^U_B}+a_1x^*+b_1=\tpi(x^*), \vspace{1.5mm} \nonumber \\
\beta^L_AA^L_1(x^*)^{\beta^L_A-1}+\beta^L_BB^L_1(x^*)^{\beta^L_B-1}=\beta^U_AA^U_1(x^*)^{\beta^U_A-1}+\beta^U_BB^U_1(x^*)^{\beta^U_B-1}+a_1, \vspace{1.5mm} \nonumber \\
\beta^L_A(\beta^L_A-1)A^L_1(x^*)^{\beta^L_A-2}+\beta^L_B(\beta^L_B-1)B^L_1(x^*)^{\beta^L_B-2} \nonumber \\
\hspace{25mm} =\beta^U_A(\beta^U_A-1)A^U_1(x^*)^{\beta^U_A-2}+\beta^U_B(\beta^U_B-1)B^U_1(x^*)^{\beta^U_B-2}. \nonumber
\end{numcases}
Solving the above, together with (\ref{eq-ABL}) and (\ref{eq-ABU}), we obtain (\ref{eq-AB1}).

%%%%%%%%%%%%%%%%%%%%%%%%%%%%%%%%%%%%%%%%%%%%%%%%%%%%%%%%%%%%%%%%%%%%%%%%%%%%%%%
\noindent{\bf Step 4:} \ In this step, we shall derive (\ref{eq-x*}). Since $V_0$ and $V^\prime_0$ are continuous at $x^*$, we have
\begin{numcases}{}
A^L_0(x^*)^{\beta^L_A}+B^L_0(x^*)^{\beta^L_B}=A^U_0(x^*)^{\beta^U_A}+B^U_0(x^*)^{\beta^U_B}+a_0x^*+b_0. \nonumber \\
\beta^L_AA^L_0(x^*)^{\beta^L_A-1}+\beta^L_BB^L_0(x^*)^{\beta^L_B-1}=\beta^U_AA^U_0(x^*)^{\beta^U_A-1}+\beta^U_BB^U_0(x^*)^{\beta^U_B-1}+a_0, \nonumber
\end{numcases}
Using (\ref{eq-AB1}) and cancelling $a_0$, we obtain
\begin{align*}
&\l(\frac{(1-\beta^L_A)P^L_A}{G^L_0(\beta^L_A)}+\frac{(1-\beta^L_B)P^L_B}{G^L_0(\beta^L_B)}+\frac{(\beta^U_A-1)P^U_A}{G^U_0(\beta^U_A)}
      +\frac{(\beta^U_B-1)P^U_B}{G^U_0(\beta^U_B)}\r)x^* \\
&\hspace{6mm} +\frac{(1-\beta^L_A)Q^L_A}{G^L_0(\beta^L_A)}+\frac{(1-\beta^L_B)Q^L_B}{G^L_0(\beta^L_B)}+\frac{(\beta^U_A-1)Q^U_A}{G^U_0(\beta^U_A)}
      +\frac{(\beta^U_B-1)Q^U_B}{G^U_0(\beta^U_B)}+\frac{b_0}{\lambda_0}=0,
\end{align*}
and denote this as $\calP x^*+\calQ=0$.
Recall that $\beta^L_A>\zeta^{L,+}_0>\beta^L_B>1$ and $\beta^U_B<\zeta^{U,-}_0<\beta^U_A<0$.
Thus, $G^L_0(\beta^L_A),G^U_0(\beta^U_B)>0$, and $G^L_0(\beta^L_B),G^U_0(\beta^U_A)<0$ hold.
Moreover, we can see easily that $P^L_A,P^U_B<0$, $P^L_B,P^U_A>0$, $Q^L_A,Q^U_B>0$, and $Q^L_B,Q^U_A<0$.
Thus, all the terms in $\calP$ are positive, and $\calQ$ are negative.
We have then $x^*=-\displaystyle{\frac{\calQ}{\calP}}>0$, that is, (\ref{eq-x*}) holds.

%%%%%%%%%%%%%%%%%%%%%%%%%%%%%%%%%%%%%%%%%%%%%%%%%%%%%%%%%%%%%%%%%%%%%%%%%%%%%%%
\noindent{\bf Step 5:} \ We show that $V_i, i=0,1$ are $\bbR_+$-valued in this last step.
Since $V_i(x)\sim a_ix+b_i$ as $x\to\infty$ and $a_i>0$ for $i=0,1$, there is an $M>0$ such that $V_i(x)>0$ for any $x>M$ and $i=0,1$.
Now, we denote
\[
V_{\oli}(\olx):=\min_{x\in(0,M]}\min_{i=0,1}V_i(x)
\]
and assume that $V_{\oli}(\olx)<0$. We have then $V^\prime_{\oli}(\olx)=0$, $V^{\prime\prime}_{\oli}(\olx)>0$, and $V_{\oli}(\olx)\leq V_{1-\oli}(\olx)$.
When $\olx\in(0,x^*)$, it follows that
\begin{equation}\label{eq-ODE-Step5}
-rV_{\oli}(\olx)+\mu_{\oli}\olx V_{\oli}^\prime(\olx)+\frac{1}{2}\sigma_{\oli}^2\olx^2V_{\oli}^{\prime\prime}(\olx)+\lambda_{\oli}(V_{1-\oli}(\olx)-V_{\oli}(\olx))= 0.
\end{equation}
Thus, we have $V_{\oli}(\olx)\geq0$, which contradicts to the assumption that $V_{\oli}(\olx)<0$. 
Next, consider the case where $\olx>x^*$. If $\oli=0$, then (\ref{eq-ODE-Step5}) holds. This is a contradiction.
When $\oli=1$, we have
\[
-rV_1(\olx)+\mu_1\olx V_1^\prime(\olx)+\frac{1}{2}\sigma_1^2\olx^2V_1^{\prime\prime}(\olx)+\lambda_1(V_0(\olx)-V_1(\olx))+\eta(\tpi(\olx)-V_1(\olx))= 0.
\]
The second, third and fourth terms are non-negative. In addition, the fifth term is also non-negative since $\tpi(\olx)>\tpi(x^*)=V_1(x^*)\geq V_1(\olx)$.
Thus, $V_1(\olx)\geq0$, which is a contradiction.
Lastly, when $\olx=x^*$, for any $\ve>0$, there is a $\delta>0$ such that $\mu_{\oli}V^\prime_{\oli}(x)>-\ve$, $V^{\prime\prime}_{\oli}(x)>0$
and $V_{1-\oli}(x)-V_{\oli}(x)>-\ve$ hold for any $x\in(x^*-\delta,x^*)$.
We have then $-rV_{\oli}(x)-\ve x^*-\lambda_{\oli}\ve\leq0$ for any $x\in(x^*-\delta,x^*)$ from the view of (\ref{eq-ODE-Step5}), which means that $V_{\oli}(x^*)\geq0$ holds.
This is a contradiction.
Consequently, $V_i, i=0,1$ are $\bbR_+$-valued. In particular, we have $V_1(x^*)=\tpi(x^*)\geq0$, from which $x^*\geq\tK$ follows.
Thus,  $V_i, i=0,1$ satisfy (\ref{VI-ODE1-2}) and (\ref{VI-*}) since $\tK\geq K$ and $\tpi(x)=\pi(x)$ for any $x\geq K$.
Consequently, $(V_0,V_1,x^*)$ gives the unique solution to Problem \ref{prob-preVI}.
This completes the proof of Proposition \ref{prop-preVI}.
\fin

%%%%%%%%%%%%%%%%%%%%%%%%%%%%%%%%%%%%%%%%%%%%%%%%%%%%%%%%%%%%%%%%%%%%%%%%%%%%%%%
\begin{rem}\label{rem3-1}
It is very complicated to show that the function $V_1$ satisfies the boundary conditions (\ref{VI-*L}) and (\ref{VI-*U}).
However, we can confirm that the conditions are met by implementing numerical computation for many parameter sets.
In fact, with $r=0.1$ and $\pi(x)=(x-0.9)^+-0.1$ fixed, and the values of $\mu_0$ and $\mu_1$ as $-10,-5, -2, -1, -0.5, 0, 0.05, 0.099$,
and $\sigma_0$, $\sigma_1$, $\lambda_0$, $\lambda_1$ and $\eta$ as $0.1, 1, 2,5$, $V_1$ satisfies (\ref{VI-*L}) and (\ref{VI-*U}) for all $65536$ parameter sets.
Thus, we can expect the boundary conditions (\ref{VI-*L}) and (\ref{VI-*U}) to be satisfied for any parameter set. We leave making sure of this fact to future research.
\end{rem}

%%%%%%%%%%%%%%%%%%%%%%%%%%%%%%%%%%%%%%%%%%%%%%%%%%%%%%%%%%%%%%%%%%%%%%%%%%%%%%%
%
% Section 4
%
%%%%%%%%%%%%%%%%%%%%%%%%%%%%%%%%%%%%%%%%%%%%%%%%%%%%%%%%%%%%%%%%%%%%%%%%%%%%%%%
\setcounter{equation}{0}
\section{Verification}
In this section, we show that the functions $V_i, i=0,1$ given in Proposition \ref{prop-preVI} coincide with the value functions $v_i, i=0,1$ defined by (\ref{def-vi}),
and an optimal stopping time $\tau^*$ exists as a threshold type with the optimal threshold $x^*$ given in (\ref{eq-x*}).
To this end, we assume that $V_1$ satisfies the boundary conditions (\ref{VI-*L}) and (\ref{VI-*U}).

Let us start with some preparations. First of all,  it is immediately apparent that the following lemma holds.

%%%%%%%%%%%%%%%%%%%%%%%%%%%%%%%%%%%%%%%%%%%%%%%%%%%%%%%%%%%%%%%%%%%%%%%%%%%%%%%
\begin{lem}\label{lem-V}
For $i=0,1$, $V^\prime_i$ is bounded, and there is a $c_i>0$ such that $V_i(x)\leq c_ix$ for any $x>0$.
\end{lem}

\noindent
In addition, we define
\[
T^1_k:=\inf\{t>T^1_{k-1}|\theta_t=1\mbox{ and }t=T_j\mbox{ for some }j\in\bbN\}
\]
for $k\in\bbN$ with the conventions $T^1_0\equiv0$ and $T^1_\infty\equiv\infty$.
Note that $T^1_k\in\calT$ represents the $k$th time when stopping is feasible, and $\calT$ is described as
\[
\calT=\{\tau\in\calT_0 \ | \ \mbox{ for each }\omega\in\Omega, \tau(\omega)=T^1_j(\omega)\mbox{ for some }j\in\bbN_\infty\}.
\]
Now, we define
\[
N^*:=\inf\{n\in\bbN|X_{T^1_n}\geq x^*\},
\]
with the convention $\inf\emptyset=\infty$. Note that $N^*$ is an $\bbN_\infty$-valued stopping time, where $\bbN_\infty:=\bbN\cup\{\infty\}$.
Hereafter, we write $Z\sim\exp(\lambda)$ when a random variable $Z$ follows the exponential distribution with parameter $\lambda>0$.

The following theorem is our main result.

%%%%%%%%%%%%%%%%%%%%%%%%%%%%%%%%%%%%%%%%%%%%%%%%%%%%%%%%%%%%%%%%%%%%%%%%%%%%%%%
\begin{thm}\label{thm-verif}
Suppose that $V_1$ satisfies (\ref{VI-*L}) and (\ref{VI-*U}). Then $v_i(x)=V_i(x)$ holds for any $x>0$ and $i=0,1$,
and the stopping time $\tau^*:=T^1_{N^*}\in\calT$ is optimal for the optimal stopping problem defined by (\ref{def-vi}).
\end{thm}

\proof We show this theorem by dividing five steps. \\
%%%%%%%%%%%%%%%%%%%%%%%%%%%%%%%%%%%%%%%%%%%%%%%%%%%%%%%%%%%%%%%%%%%%%%%%%%%%%%%
\noindent{\bf Step 1:} \ In this step, we fix $\theta_0=0$ and $X_0=x$, and denote $\xi_0:=\inf\{t>0|\theta_t=1\}$.
For $i=0,1$, we denote by $Y^i=\{Y^i_t\}_{t\geq0}$ a geometric Brownian motion starting at $1$ under regime $i$, that is, the solution to the following SDE:
\[
dY^i_t=Y^i_t(\mu_idt+\sigma_idW_t), \ \ \ Y^i_0=1.
\]
In the following, when we write $Y^i_t$, its independent copy may be taken if necessary.
Note that $xY^0_t=X_t$ holds if $t<\xi_0$, and $\xi_0\sim\exp(\lambda_0)$.
Now, we see the following:
\begin{equation}\label{eq-V0}
V_0(x)=\bbE\l[\int_0^\infty e^{-(r+\lambda_0)t}\lambda_0V_1(xY^0_t)dt\r], \ \ \ x>0.
\end{equation}
To this end, we define firstly
\begin{equation}\label{eq-Phi0}
\Phi^0_t:=e^{-(r+\lambda_0)t}V_0(xY^0_t){\bf 1}_{\{t<\xi_0\}}.
\end{equation}
Ito's formula implies that
\begin{align*}
\Phi^0_t &= \Bigg\{V_0(x)+\int_0^te^{-(r+\lambda_0)s}\l(-(r+\lambda_0)V_0(xY^0_s)+\calA_0V_0(xY^0_s)\r)ds \\
         & \hspace{6mm}+\int_0^te^{-(r+\lambda_0)s}\sigma_0xY^0_sV_0^\prime(xY^0_s)dW_s\Bigg\}{\bf 1}_{\{t<\xi_0\}}.
\end{align*}
Taking expectation on both sides, we have
\begin{align*}
\bbE[\Phi^0_t] &= V_0(x)\bbE[{\bf 1}_{\{t<\xi_0\}}]+\bbE\l[\int_0^te^{-(r+\lambda_0)s}(-\lambda_0)V_1(xY^0_s)ds{\bf 1}_{\{t<\xi_0\}}\r] \\
               & \hspace{6mm}+\bbE\l[\int_0^te^{-(r+\lambda_0)s}\sigma_0xY^0_sV_0^\prime(xY^0_s)dW_s{\bf 1}_{\{t<\xi_0\}}\r] \\
               &= \Bigg\{V_0(x)-\bbE\l[\int_0^te^{-(r+\lambda_0)s}\lambda_0V_1(xY^0_s)ds\r] \\
               & \hspace{6mm}+\bbE\l[\int_0^te^{-(r+\lambda_0)s}\sigma_0xY^0_sV_0^\prime(xY^0_s)dW_s\r]\Bigg\}e^{-\lambda_0t} \\
               &= \l\{V_0(x)-\bbE\l[\int_0^te^{-(r+\lambda_0)s}\lambda_0V_1(xY^0_s)ds\r]\r\}e^{-\lambda_0t}.
\end{align*}
The first equality is due to (\ref{VI-ODE0}); the second is due to the independence of $\xi_0$ and $W$ and $\xi_0\sim\exp(\lambda_0)$.
The last equality is obtained from the boundedness of $V_0^\prime$ by Lemma \ref{lem-V} and the integrability of $\displaystyle{\int_0^t(Y^0_s)^2ds}$.
From the view of (\ref{eq-Phi0}), we obtain
\[
V_0(x)=\bbE\l[e^{-(r+\lambda_0)t}V_0(xY^0_t)\r]+\bbE\l[\int_0^te^{-(r+\lambda_0)s}\lambda_0V_1(xY^0_s)ds\r].
\]
Since $V_0(x)\leq c_0x$ from Lemma \ref{lem-V} and $r>\mu_0$ from (\ref{r-mu}), we have
\[
\bbE\l[e^{-(r+\lambda_0)t}V_0(xY^0_t)\r]\leq\bbE\l[e^{-(r+\lambda_0)t}c_0xY^0_t\r],
\]
which tends to 0 as $t\to\infty$. As a result, since $V_1\geq0$, the monotone convergence theorem implies (\ref{eq-V0}).

Since $\xi_0\sim\exp(\lambda_0)$, (\ref{eq-V0}) can be rewritten as
\begin{equation}\label{eq-V0-1}
V_0(x)=\bbE\l[e^{-r\xi_0}V_1(xY^0_{\xi_0})\r]=\bbE^{0,x}\l[e^{-r\xi_0}V_1(X_{\xi_0})\r].
\end{equation}
From the view of (\ref{def-vi}), showing $v_1=V_1$, we obtain $v_0=V_0$ immediately.
In what follows, we focus on the proof of $v_1=V_1$.

%%%%%%%%%%%%%%%%%%%%%%%%%%%%%%%%%%%%%%%%%%%%%%%%%%%%%%%%%%%%%%%%%%%%%%%%%%%%%%%
\noindent{\bf Step 2:} \ Throughout the rest of this proof, we fix $\theta_0=1$ and $X_0=x$. Now, we define
\[
\olV(x):=\pi(x)\vee V_1(x)=\l\{\begin{array}{ll}\pi(x), & x\geq x^*, \\ V_1(x), & x<x^*.\end{array}\r.
\]
We can then unify (\ref{VI-ODE1-1}) and (\ref{VI-ODE1-2}) into
\begin{equation}\label{ODE-V1-2}
-rV_1(x)+\calA_1 V_1(x)+\lambda_1(V_0(x)-V_1(x))+\eta(\olV(x)-V_1(x))=0, \ \ \ x>0.
\end{equation}

Here we aim to show the following by a similar argument to Step 1:
\begin{equation}\label{eq-V1}
V_1(x)=\bbE\l[\int_0^\infty e^{-(r+\lambda_1+\eta)t}\{\lambda_1V_0(xY^1_t)+\eta\olV(xY^1_t)\}dt\r].
\end{equation}
To this end, we define
\[
\Phi^1_t:=e^{-(r+\lambda_1+\eta)t}V_1(xY^1_t){\bf 1}_{\{t<\xi_1\wedge T_1\}},
\]
where $\xi_1:=\inf\{t>0|\theta_t=0\}$. In addition, recall that $T_1=\inf\{t>0|J_t=1\}$, that is, the first investment opportunity time.
Noting that $\bbP(t<\xi_1\wedge T_1)=e^{-(\lambda_1+\eta)t}$, we obtain
\[
\bbE[\Phi^1_t]=\l\{V_1(x)-\bbE\l[\int_0^te^{-(r+\lambda_1+\eta)s}\l(\lambda_1V_0(xY^1_s)+\eta\olV(xY^1_s)\r)ds\r]\r\}e^{-(\lambda_1+\eta)t}
\]
from Ito's formula and (\ref{ODE-V1-2}). By the same sort of argument as Step 1, (\ref{eq-V1}) follows.

%%%%%%%%%%%%%%%%%%%%%%%%%%%%%%%%%%%%%%%%%%%%%%%%%%%%%%%%%%%%%%%%%%%%%%%%%%%%%%%
\noindent{\bf Step 3:} \ This step is devoted to preparing some notations.
First of all, we define two sequences of stopping times inductively as follows: $\xi^{0\to1}_0\equiv0$ and, for $k\in\bbN$, 
\begin{align*}
\xi^{1\to0}_k&:= \inf\{t>\xi^{0\to1}_{k-1}|\theta_{t-}=1,\theta_t=0\}, \\
\xi^{0\to1}_k&:= \inf\{t>\xi^{1\to0}_k|\theta_{t-}=0,\theta_t=1\}.
\end{align*}
We call the time interval $[\xi^{0\to1}_{k-1},\xi^{0\to1}_k)$ the $k$th phase.
Note that each phase begins when the regime changes into $1$, moves to regime $0$ midway through, and ends when it returns to regime $1$ again.
Moreover, we define the following two sequences of i.i.d. random variables:
\[
U^1_k:= \xi^{1\to0}_k-\xi^{0\to1}_{k-1}, \ \ \ U^0_k:= \xi^{0\to1}_k-\xi^{1\to0}_k.
\]
Note that each $U^i_k\sim\exp(\lambda_i)$ expresses the length of regime $i$ in the $k$th phase, and $U^0_{k_0}$ and $U^1_{k_1}$ are independent for any $k_0,k_1\in\bbN$.
For $k\in\bbN$, we denote by $\tT_k$ the first investment opportunity time after the start of the $k$th phase, that is,
\[
\tT_k:=\inf\{t>\xi^{0\to1}_{k-1}|t=T_j\mbox{ for some }j\in\bbN\}.
\]
Note that $\tT_k$ is not necessarily in the $k$th phase, and $\theta_{\tT_k}$ may take the value of $0$.
In addition, we define $U^P_k:=\tT_k-\xi^{0\to1}_{k-1}\sim\exp(\eta)$,
which represents the length of time from the start of the $k$th phase until the arrival of the first investment opportunity.

%%%%%%%%%%%%%%%%%%%%%%%%%%%%%%%%%%%%%%%%%%%%%%%%%%%%%%%%%%%%%%%%%%%%%%%%%%%%%%%
\noindent{\bf Step 4:} \ In this step, we shall show 
\begin{equation}\label{eq-Step4}
V_1(x)=\bbE^{1,x}\l[e^{-rT^1_1}\olV(X_{T^1_1})\r].
\end{equation}
Recall $T^1_1=\inf\{t>0|\theta_t=1\mbox{ and }t=T_j\mbox{ for some }j\in\bbN\}$, that is, the time when stopping becomes feasible for the first time.
 
First of all, we can rewrite (\ref{eq-V1}) as
\begin{equation}\label{eq-V1-1}
V_1(x)=\bbE\l[e^{-rU^1_1}V_0(xY^1_{U^1_1}){\bf 1}_{\{U^1_1<U^P_1\}}+e^{-rU^P_1}\olV(xY^1_{U^P_1}){\bf 1}_{\{U^P_1<U^1_1\}}\r],
\end{equation}
since $U^1_1$ is independent of $U^P_1$ and $\bbP(U^P_1>t)=e^{-\eta t}$. Using (\ref{eq-V0-1}) and (\ref{eq-V1-1}), we have
\begin{align*}
V_1(x)
&= \bbE\l[e^{-rU^1_1}\l(e^{-rU^0_1}V_1(xY^1_{U^1_1}Y^0_{U^0_1})\r){\bf 1}_{\{U^1_1<U^P_1\}}+e^{-rU^P_1}\olV(xY^1_{U^P_1}){\bf 1}_{\{U^P_1<U^1_1\}}\r] \\
&= \bbE\Big[e^{-r(U^1_1+U^0_1)}\Big(e^{-rU^1_2}V_0(xY^1_{U^1_1}Y^0_{U^0_1}Y^1_{U^1_2}){\bf 1}_{\{U^1_2<U^P_2\}} \\
&  \hspace{6mm}+e^{-rU^P_2}\olV(xY^1_{U^1_1}Y^0_{U^0_1}Y^1_{U^P_2}){\bf 1}_{\{U^P_2<U^1_2\}}\Big){\bf 1}_{\{U^1_1<U^P_1\}}+e^{-rU^P_1}\olV(xY^1_{U^P_1}){\bf 1}_{\{U^P_1<U^1_1\}}\Big] \\
&= \bbE\Big[e^{-r(U^1_1+U^0_1+U^1_2)}V_0(xY^1_{U^1_1}Y^0_{U^0_1}Y^1_{U^1_2}){\bf 1}_{\{U^1_1<U^P_1\}\cap\{U^1_2<U^P_2\}} \\
&  \hspace{6mm}+e^{-r(U^1_1+U^0_1+U^P_2)}\olV(xY^1_{U^1_1}Y^0_{U^0_1}Y^1_{U^P_2}){\bf 1}_{\{U^1_1<U^P_1\}\cap\{U^P_2<U^1_2\}}
   +e^{-rU^P_1}\olV(xY^1_{U^P_1}){\bf 1}_{\{U^P_1<U^1_1\}}\Big].
\end{align*}
Note that all random variables in the above are independent.
Now, we denote
\[
Z^0_n:=\exp\l\{-r\l(\sum_{k=1}^nU^1_k+\sum_{k=1}^{n-1}U^0_k\r)\r\}V_0\l(x\prod_{k=1}^{n-1}\l(Y^1_{U^1_k}Y^0_{U^0_k}\r)Y^1_{U^1_n}\r){\bf 1}_{\bigcap_{k=1}^n\{U^1_k<U^P_k\}}
\]
for $n\in\bbN$, $\olZ_1:=\displaystyle{e^{-rU^P_1}\olV(xY^1_{U^P_1}){\bf 1}_{\{U^P_1<U^1_1\}}}$, and
\[
\olZ_k:=\exp\l\{-r\l(\sum_{j=1}^{k-1}(U^1_j+U^0_j)+U^P_k\r)\r\}\olV\l(x\prod_{j=1}^{k-1}\l(Y^1_{U^1_j}Y^0_{U^0_j}\r)Y^1_{U^P_k}\r)
        {\bf 1}_{\bigcap_{j=1}^{k-1}\{U^1_j<U^P_j\}\cap\{U^P_k<U^1_k\}}
\]
for $k\geq2$.
Remark that, for $k\in\bbN$, we can rewrite $\olZ_k$ as follows:
\begin{equation}\label{eq-olZ}
\olZ_k=e^{-rT^1_1}\olV(X_{T^1_1}){\bf 1}_{\{\xi^{0\to1}_{k-1}\leq T^1_1<\xi^{1\to0}_k\}}
\end{equation}
when $\theta_0=1$ and $X_0=x$. We have then, for any $n\in\bbN$,
\[
V_1(x)=\bbE\l[Z^0_n+\sum_{k=1}^n\olZ_k\r].
\]

From Lemma \ref{lem-V} and the independence of all random variables, it follows that
\begin{align*}
\bbE[Z^0_n]
&\leq \bbE\l[\exp\l\{-r\l(\sum_{k=1}^nU^1_k+\sum_{k=1}^{n-1}U^0_k\r)\r\}V_0\l(x\prod_{k=1}^nY^1_{U^1_k}\prod_{k=1}^{n-1}Y^0_{U^0_k}\r)\r] \\
&\leq \bbE\l[c_0x\prod_{k=1}^n\l(e^{-rU^1_k}Y^1_{U^1_k}\r)\prod_{k=1}^{n-1}\l(e^{-rU^0_k}Y^1_{U^0_k}\r)\r]
      =c_0x\prod_{k=1}^n\bbE\l[e^{-rU^1_k}Y^1_{U^1_k}\r]\prod_{k=1}^{n-1}\bbE\l[e^{-rU^0_k}Y^0_{U^0_k}\r] \\
&\leq c_0x\l(\frac{\lambda_1}{r-\mu_1+\lambda_1}\r)^n\l(\frac{\lambda_0}{r-\mu_0+\lambda_0}\r)^{n-1}
\end{align*}
since
\[
\bbE\l[e^{-rU^i_k}Y^i_{U^i_k}\r]=\frac{\lambda_i}{r-\mu_i+\lambda_i}.
\]
As a result, we obtain $\displaystyle{\lim_{n\to\infty}\bbE[Z^0_n]=0}$.
Since each $\olZ_k$ is non-negative, the monotone convergence theorem implies that
\[
V_1(x)=\lim_{n\to\infty}\bbE\l[Z^0_n+\sum_{k=1}^n\olZ_k\r]=\bbE\l[\sum_{k=1}^\infty\olZ_k\r].
\]
Thus, (\ref{eq-olZ}) provides that
\[
V_1(x)=\bbE^{1,x}\l[e^{-rT^1_1}\olV(X_{T^1_1})\sum_{k=1}^\infty{\bf 1}_{\{\xi^{0\to1}_{k-1}\leq T^1_1<\xi^{1\to0}_k\}}\r]=\bbE^{1,x}\l[e^{-rT^1_1}\olV(X_{T^1_1}){\bf 1}_{\{T^1_1<\infty\}}\r].
\]
On the other hand, $e^{-rt}\olV(X_t)\leq e^{-rt}(c_1\vee\alpha)X_t$ holds.
Since $e^{-rt}X_t$ is a non-negative supermartingale, it converges to $0$ a.s. as $t\to\infty$ by, e.g., Problem 1.3.16 of \cite{KS12}. As a result, we have
\[
\bbE^{1,x}\l[e^{-rT^1_1}\olV(X_{T^1_1}){\bf 1}_{\{T^1_1<\infty\}}\r]=\bbE^{1,x}\l[e^{-rT^1_1}\olV(X_{T^1_1})\r],
\]
from which (\ref{eq-Step4}) follows.

%%%%%%%%%%%%%%%%%%%%%%%%%%%%%%%%%%%%%%%%%%%%%%%%%%%%%%%%%%%%%%%%%%%%%%%%%%%%%%%
\noindent{\bf Step 5:} \ We define a filtration $\bbG=\{\calG_n\}_{n\in\bbN_0}$ as $\calG_n:=\calF_{T^1_n}$
and  a process $\olS=\{\olS_n\}_{n\in\bbN_0}$ as $\olS_n:=e^{-rT^1_n}\olV(X_{T^1_n})$, where $\bbN_0:=\bbN\cup\{0\}$.
We have then, for any $n\in\bbN_0$,
\begin{align*}
\olS_n &\geq e^{-rT^1_n}V_1(X_{T^1_n})=e^{-rT^1_n}\bbE^{1,y}\l[e^{-r\hT^1_1}\olV(X_{\hT^1_1})\r]\Bigg|_{y=X_{T^1_n}} \\
      &= \bbE^{1,x}\l[e^{-rT^1_{n+1}}\olV(X_{T^1_{n+1}})\Big|\calG_n\r]=\bbE^{1,x}\l[\olS_{n+1}|\calG_n\r],
\end{align*}
where $\hT^1_1$ is an independent copy of $T^1_1$.
Thus, $\olS$ is a non-negative $\bbG$-supermartingale, and $\olS_n$ converges to $0$ a.s. as $n\to\infty$.
On the other hand, Lemma 1 of \cite{DW02} implies
\begin{equation}\label{eq-step5-1}
\calT=\{T^1_N|N\mbox{ is an $\bbN_\infty$-valued $\bbG$-stopping time}\}.
\end{equation}
Since $\olS_0\geq\bbE^{1,x}\l[\olS_n\r]\geq0$ for any $n\in\bbN$, the optional sampling theorem, e.g., Theorem 16 of Chapter V in \cite{DM82}, together with (\ref{eq-Step4}), yields that
\[
V_1(x)=\bbE^{1,x}\l[\olS_1\r]\geq\bbE^{1,x}\l[\olS_N\r]\geq\bbE^{1,x}\l[e^{-rT^1_N}\pi(X_{T^1_N})\r]
\]
for any $\bbN_\infty$-valued $\bbG$-stopping time $N$.
Taking supremum on the right-hand side over all such $N$'s, we obtain $v_1\leq V_1$ from the view of (\ref{def-vi}) and (\ref{eq-step5-1}).

Next, we see the reverse inequality $v_1\geq V_1$.
To this end, we recall $N^*:=\inf\{n\in\bbN|X_{T^1_n}\geq x^*\}$ and define $\olS^*_n:=\exp\{-rT^1_{N^*\wedge n}\}\olV(X_{T^1_{N^*\wedge n}})$ for $n\in\bbN_0$.
As shown in Lemma \ref{lem-olS*}, $\olS^*=\{\olS^*_n\}_{n\in\bbN_0}$ is a uniformly integrable martingale, which implies that
\begin{align*}
V_1(x) &\leq \olV(x)=\olS^*_0=\lim_{n\to\infty}\bbE^{1,x}\l[\olS^*_n\r]=\bbE^{1,x}\l[\lim_{n\to\infty}\olS^*_n\r] \\
       &= \bbE^{1,x}\l[e^{-rT^1_{N^*}}\olV(X_{T^1_{N^*}})\r]=\bbE^{1,x}\l[e^{-rT^1_{N^*}}\pi(X_{T^1_{N^*}})\r]\leq v_1(x)
\end{align*}
since $\olV(x)=\pi(x)$ for any $x\geq x^*$, and $T^1_{N^*}\in\calT$. Consequently, we obtain
\[
v_1(x)=V_1(x)=\bbE^{1,x}\l[e^{-rT^1_{N^*}}\pi(X_{T^1_{N^*}})\r], \ \ \ x>0,
\]
and thus, the stopping time $T^1_{N^*}$ is optimal. This completes the proof of Theorem \ref{thm-verif}.
\fin

%%%%%%%%%%%%%%%%%%%%%%%%%%%%%%%%%%%%%%%%%%%%%%%%%%%%%%%%%%%%%%%%%%%%%%%%%%%%%%%
\begin{lem}\label{lem-olS*}
$\olS^*$ is a uniformly integrable martingale.
\end{lem}

\proof
We shall prove this lemma by the same sort of argument as Step 2 of Section 3.2 in \cite{DW02}.
First of all, for any $n\in\bbN$, we have
\begin{align*}
\bbE^{1,x}\l[\olS^*_n\big|\calG_{n-1}\r]
&= \bbE^{1,x}\l[e^{-rT^1_n}\olV(X_{T^1_n}){\bf 1}_{\{N^*\geq n\}}\big|\calG_{n-1}\r]+\bbE^{1,x}\l[e^{-rT^1_{N^*}}\olV(X_{T^1_{N^*}}){\bf 1}_{\{N^*<n\}}\big|\calG_{n-1}\r] \\
&= e^{-rT^1_{n-1}}\bbE^{1,y}\l[e^{-r\hT^1_1}\olV(X_{\hT^1_1})\r]\Bigg|_{y=X_{T^1_{n-1}}}\hspace{-5mm}{\bf 1}_{\{N^*\geq n\}}+e^{-rT^1_{N^*}}\olV(X_{T^1_{N^*}}){\bf 1}_{\{N^*<n\}} \\
&= e^{-rT^1_{n-1}}V_1(X_{T^1_{n-1}}){\bf 1}_{\{N^*\geq n\}}+e^{-rT^1_{N^*}}\olV(X_{T^1_{N^*}}){\bf 1}_{\{N^*<n\}} \\
&= e^{-rT^1_{n-1}}\olV(X_{T^1_{n-1}}){\bf 1}_{\{N^*\geq n\}}+e^{-rT^1_{N^*}}\olV(X_{T^1_{N^*}}){\bf 1}_{\{N^*<n\}}=\olS^*_{n-1},
\end{align*}
where $\hT^1_1$ is an independent copy of $T^1_1$. As a result, $\olS^*$ is a $\bbG$-martingale.

Next, we show the uniform integrability. To see this, we have only to show that
\[
\sup_{n\in\bbN}\bbE^{1,x}\l[\l|\olS^*_n\r|^p\r]<\infty \ \ \ \mbox{for some }p>1.
\]
Since $\olV(x)\leq(c_1\vee\alpha)x$, it suffices to see that
\begin{equation}\label{eq-lem-olS*}
\sup_{n\in\bbN}\bbE^{1,x}\l[\exp\l\{-prT^1_{N^*\wedge n}\r\}X^p_{T^1_{N^*\wedge n}}\r]<\infty \ \ \ \mbox{for some }p>1.
\end{equation}
Note that
\begin{align*}
e^{-prt}X^p_t
&= x^p\exp\l\{p\int_0^t\l(\mu_{\theta_s}-r-\frac{1}{2}\sigma^2_{\theta_s}\r)ds+p\int_0^t\sigma_{\theta_s}dW_s\r\} \\
&= x^p\exp\l\{p\int_0^t\l(\mu_{\theta_s}-r+\frac{p-1}{2}\sigma^2_{\theta_s}\r)ds-\int_0^t\frac{p^2}{2}\sigma^2_{\theta_s}ds+\int_0^tp\sigma_{\theta_s}dW_s\r\}.
\end{align*}
Now, we take a $p>1$ satisfying $\mu_i-r+\frac{\sigma_i^2}{2}(p-1)<0$ for any $i=0,1$. Denoting
\[
M^*_n:=e^{-prT^1_n}X^p_{T^1_n}, \ \ \ n\in\bbN_0,
\]
we can see that $M^*=\{M^*_n\}_{n\in\bbN_0}$ is a nonnegative $\bbG$-supermartingale.
Thus, the optional sampling theorem, e.g., Theorem 16 of Chapter V in \cite{DM82}, implies that
\[
\bbE^{1,x}\l[\exp\l\{-prT^1_{N^*\wedge n}\r\}X^p_{T^1_{N^*\wedge n}}\r]=\bbE^{1,x}\l[M^*_{N^*\wedge n}\r]\leq M^*_0=x^p
\]
holds for any $n\in\bbN$, from which (\ref{eq-lem-olS*}) follows.
\fin

\begin{figure}[H]
\begin{minipage}{0.5\hsize}
    \includegraphics[width=70mm]{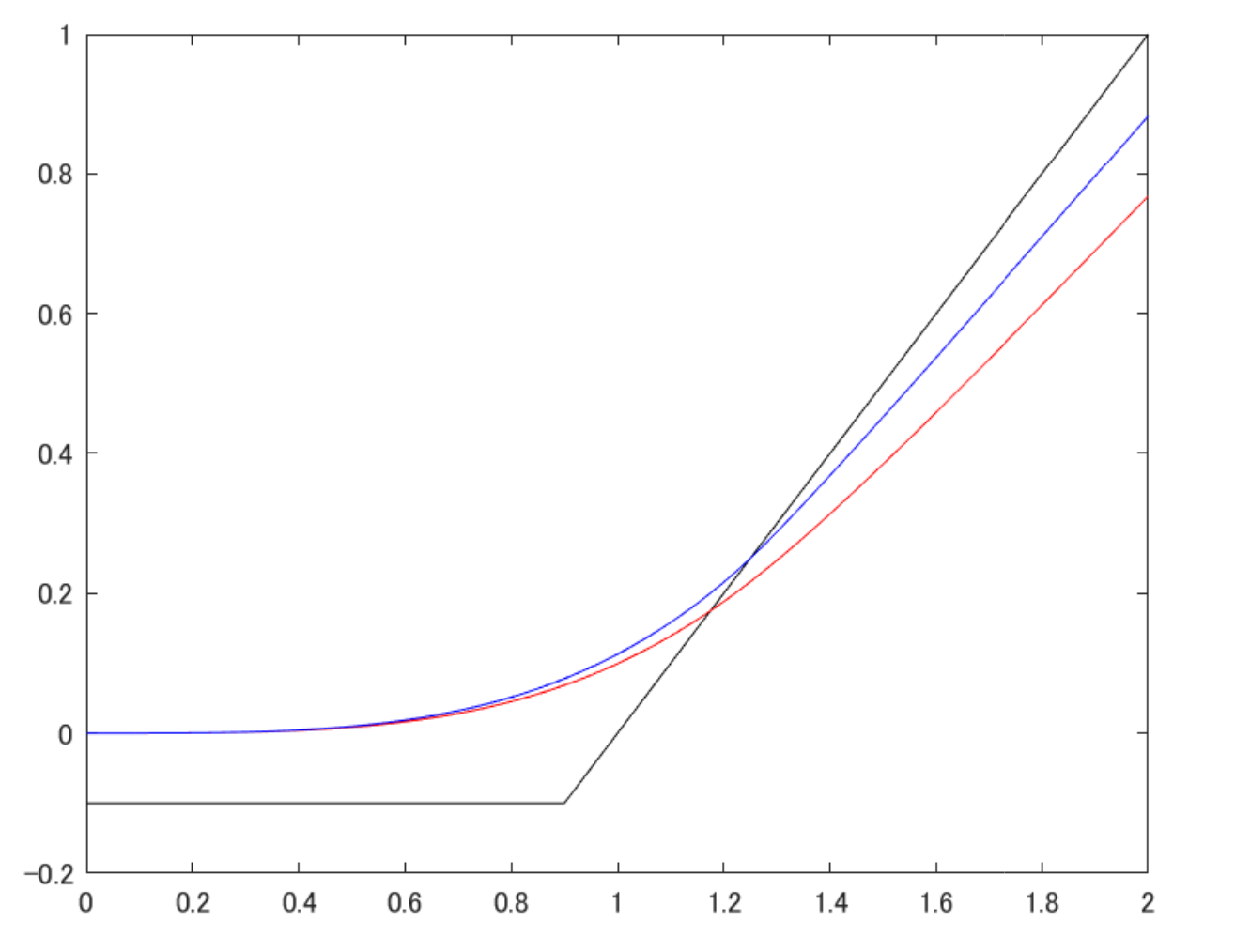}\vspace{-2mm}\caption{}\label{fig1}
\end{minipage}
\begin{minipage}{0.5\hsize}
    \includegraphics[width=70mm]{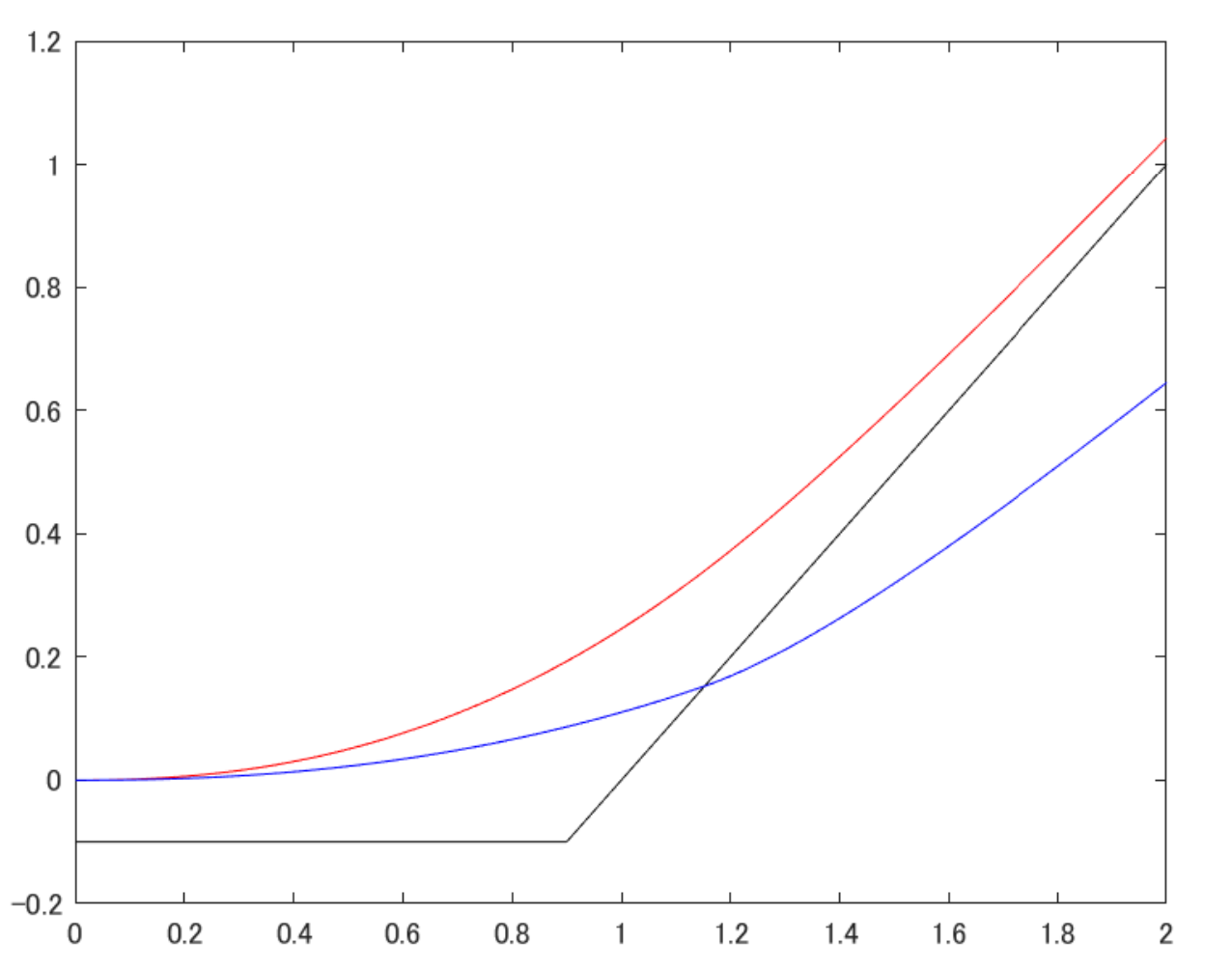}\vspace{-2mm}\caption{}\label{fig2}
\end{minipage}
\end{figure}

By Theorem \ref{thm-verif}, an optimal stopping time $\tau^*$ exists as a threshold type with the optimal threshold $x^*$
if $V_1$ in Proposition \ref{prop-preVI} satisfies the boundary conditions (\ref{VI-*L}) and (\ref{VI-*U}).
Moreover, (\ref{eq-ViL}), (\ref{eq-ViU}) and (\ref{eq-x*}) give expressions of the value functions $v_i, i=0,1$ and the optimal threshold $x^*$, respectively.
Although these expressions contain solutions to quartic equations, we can compute the value of $x^*$ numerically and illustrate the value functions $v_i, i=0,1$, e.g.,
for the case where $\pi(x)=(x-0.9)^+-0.1$, $r=0.1$, $\mu_0=-0.1$, $\mu_1=0.05$, $\sigma_0=0.2$, $\sigma_1=0.1$, $\lambda_0=2$, $\lambda_1=1$ and $\eta=1$, we obtain approximately
\[
v_0(x)=\l\{\begin{array}{l}
-4.05\times10^{-5}x^{17.18}+0.10x^{3.52}, \ \ \ 0<x<x^*, \\
0.16x^{-5.28}-0.12x^{-26.12}+0.80x-0.83, \ \ \ x>x^*,
\end{array}\r.
\]
\[
v_1(x)=\l\{\begin{array}{l}
-3.53\times10^{-5}x^{17.18}+0.11x^{3.52}, \ \ \ 0<x<x^*, \\
0.07x^{-5.28}-0.91x^{-26.12}+0.88x-0.87, \ \ \ x>x^*,
\end{array}\r.
\]
and $x^*=1.250142442232948$.
Figure \ref{fig1} illustrates the functions $v_0(x)$, $v_1(x)$ and $\pi(x)$ by red, blue, and black curves.
Furthermore, it is immediately seen that the value functions $v_i, i=0,1$ are non-negative non-decreasing convex functions and $v_i(x)\sim a_ix$ as $x\to\infty$ for $i=0,1$.
However, the magnitude relationship of $v_0$ and $v_1$ depends on how we take parameters.
The function $v_1$ is larger in the above example but simply replacing the values of $\mu_0$ and $\mu_1$ with $0.5$ and $-0.5$, respectively, reverses the magnitude relationship
between $v_0$ and $v_1$ as illustrated in Figure \ref{fig2}. Besides, $x^*$ for this case takes the value of $1.152507688970727$.

%%%%%%%%%%%%%%%%%%%%%%%%%%%%%%%%%%%%%%%%%%%%%%%%%%%%%%%%%%%%%%%%%%%%%%%%%%%%%%%

%%%%%%%%%%%%%%%%%%%%%%%%%%%%%%%%%%%%%%%%%%%%%%%%%%%%%%%%%%%%%%%%%%%%%%%%%%%%%%%
%
% Section 5
%
%%%%%%%%%%%%%%%%%%%%%%%%%%%%%%%%%%%%%%%%%%%%%%%%%%%%%%%%%%%%%%%%%%%%%%%%%%%%%%%
\setcounter{equation}{0}
\section{Asymptotic behaviors}
This section discusses asymptotic behaviors of the value functions $v_i, i=0,1$ and the optimal threshold $x^*$ when some parameter goes to $\infty$.
To compare with results in preceding literature, we consider the case where $X$ is a geometric Brownian motion given as $dX_t=X_t(\mu dt+\sigma dW_t)$, that is,
$\mu=\mu_0=\mu_1$ and $\sigma= \sigma_0=\sigma_1$. Then, simple calculations show that
\begin{equation}\label{eq-Sect5-1}
\l\{\begin{array}{l}
\beta^L_A=\frac{1}{2}-\frac{\mu}{\sigma^2}+\sqrt{\l(\frac{1}{2}-\frac{\mu}{\sigma^2}\r)^2+\frac{2(\lambda_0+\lambda_1+r)}{\sigma^2}}, \ \ \ 
\beta^L_B=\frac{1}{2}-\frac{\mu}{\sigma^2}+\sqrt{\l(\frac{1}{2}-\frac{\mu}{\sigma^2}\r)^2+\frac{2r}{\sigma^2}}, \vspace{2mm} \\
\beta^U_A=\frac{1}{2}-\frac{\mu}{\sigma^2}-\sqrt{\l(\frac{1}{2}-\frac{\mu}{\sigma^2}\r)^2
           +\frac{1}{\sigma^2}\l(\lambda_0+\lambda_1+\eta+2r-\sqrt{(\lambda_0+\lambda_1+\eta)^2-4\lambda_0\eta}\r)}, \vspace{2mm} \\
\beta^U_B=\frac{1}{2}-\frac{\mu}{\sigma^2}-\sqrt{\l(\frac{1}{2}-\frac{\mu}{\sigma^2}\r)^2
           +\frac{1}{\sigma^2}\l(\lambda_0+\lambda_1+\eta+2r+\sqrt{(\lambda_0+\lambda_1+\eta)^2-4\lambda_0\eta}\r)}.
\end{array}\r.
\end{equation}

%%%%%%%%%%%%%%%%%%%%%%%%%%%%%%%%%%%%%%%%%%%%%%%%%%%%%%%%%%%%%%%%%%%%%%%%%%%%%%%
\subsection{Asymptotic behaviors as $\eta\to\infty$}
When $\eta\to\infty$, investment opportunities arrive continuously, which means only the regime constraint remains.
First of all, we have
\begin{equation}\label{eq-Sect5-2}
\lim_{\eta\to\infty}\beta^U_A=\frac{1}{2}-\frac{\mu}{\sigma^2}-\sqrt{\l(\frac{1}{2}-\frac{\mu}{\sigma^2}\r)^2+\frac{2(\lambda_0+r)}{\sigma^2}}=\zeta^{L,-}_0,
\ \mbox{ and } \ \lim_{\eta\to\infty}\beta^U_B=-\infty,
\end{equation}
respectively, but the values of $\beta^L_A$ and $\beta^L_B$ are independent of $\eta$.
In addition, it follows that
\begin{equation}\label{eq-Sect5-3}
a_0\to\frac{\alpha\lambda_0}{r-\mu+\lambda_0}, \ \ \ a_1\to\alpha, \ \ \ b_0\to-\frac{\alpha\tK\lambda_0}{r+\lambda_0}, \ \mbox{ and } \ b_1\to-\alpha\tK
\end{equation}
as $\eta\to\infty$. By (\ref{eq-PQ}) and (\ref{eq-Sect5-2}), we  can see that
\[
P^L_A\to\frac{(-\beta^L_B+1)\alpha}{\beta^L_A-\beta^L_B}, \ \ \ Q^L_A\to\frac{\beta^L_B\alpha\tK}{\beta^L_A-\beta^L_B}, \ \ \ 
P^L_B\to\frac{(\beta^L_A-1)\alpha}{\beta^L_A-\beta^L_B}, \ \mbox{ and } \ Q^L_B\to\frac{-\beta^L_A\alpha\tK}{\beta^L_A-\beta^L_B}
\]
as $\eta\to\infty$, and $P^U_A,Q^U_A,P^U_B,Q^U_B$ converge to 0.
By Proposition \ref{prop-preVI}, Theorem \ref{thm-verif} and (\ref{eq-Sect5-1}), we obtain
\[
\lim_{\eta\to\infty}v_1(x)= \l\{\begin{array}{l}
          \displaystyle{\frac{(-\beta^L_B+1)\alpha x^*_\infty+\beta^L_B\alpha\tK}{\beta^L_A-\beta^L_B}\l(\frac{x}{x^*_\infty}\r)^{\beta^L_A}
             +\frac{(\beta^L_A-1)\alpha x^*_\infty-\beta^L_A\alpha\tK}{\beta^L_A-\beta^L_B}\l(\frac{x}{x^*_\infty}\r)^{\beta^L_B}}, \ \ \ 0<x<x^*_\infty, \vspace{1mm} \\
          \alpha x -\alpha\tK, \ \ \ x>x^*_\infty,
        \end{array}\r.
\]
where $\beta^L_A$ and $\beta^L_B$ are given in (\ref{eq-Sect5-1}), and $x^*_\infty:=\displaystyle{\lim_{\eta\to\infty}x^*}$ given in (\ref{eq-Sect5-5}) below.
Now, we assume that $\displaystyle{\lim_{\eta\to\infty}v_1(x)}\geq\pi(x)$ for any $x\in(0,x^*_\infty)$.
Since $G^L_0(\beta^L_A)=\lambda_1$ and $G^L_0(\beta^L_B)=-\lambda_0$, we have
\[
\lim_{\eta\to\infty}v_0(x)= -\frac{\lambda_0}{\lambda_1}\frac{(-\beta^L_B+1)\alpha x^*_\infty+\beta^L_B\alpha\tK}{\beta^L_A-\beta^L_B}\l(\frac{x}{x^*_\infty}\r)^{\beta^L_A}
             +\frac{(\beta^L_A-1)\alpha x^*_\infty-\beta^L_A\alpha\tK}{\beta^L_A-\beta^L_B}\l(\frac{x}{x^*_\infty}\r)^{\beta^L_B}, \ \ \ 0<x<x^*_\infty,
\]
by (\ref{eq-ABL}). In addition, the continuity of $V_0$ at $x^*_\infty$, togther with (\ref{eq-Sect5-3}) and $P^U_B,Q^U_B\to0$, implies that
\[
\lim_{\eta\to\infty}v_0(x)= \olA^U_0\l(\frac{x}{x^*_\infty}\r)^{\zeta^{L,-}_0}+\frac{\alpha\lambda_0}{r-\mu+\lambda_0}x-\frac{\alpha\tK\lambda_0}{r+\lambda_0}, \ \ \ x>x^*_\infty,
\]
where $\zeta^{L,-}_0=\displaystyle{\lim_{\eta\to\infty}\beta^U_A}$ by (\ref{eq-Sect5-2}), and
\begin{equation}\label{eq-Sect5-4}
\olA^U_0:= -\frac{\lambda_0}{\lambda_1}\frac{(-\beta^L_B+1)\alpha x^*_\infty+\beta^L_B\alpha\tK}{\beta^L_A-\beta^L_B}
                +\frac{(\beta^L_A-1)\alpha x^*_\infty-\beta^L_A\alpha\tK}{\beta^L_A-\beta^L_B}-\frac{\alpha\lambda_0}{r-\mu+\lambda_0}x^*_\infty+\frac{\alpha\tK\lambda_0}{r+\lambda_0}.
\end{equation}
From the view of (\ref{eq-Sect5-4}), we have
\[
\lim_{\eta\to\infty}\frac{-\lambda_0P^U_A}{G^U_0(\beta^U_A)}=-\frac{\lambda_0}{\lambda_1}\frac{(-\beta^L_B+1)\alpha}{\beta^L_A-\beta^L_B}
                +\frac{(\beta^L_A-1)\alpha}{\beta^L_A-\beta^L_B}-\frac{\alpha\lambda_0}{r-\mu+\lambda_0}
\]
and
\[
\lim_{\eta\to\infty}\frac{-\lambda_0Q^U_A}{G^U_0(\beta^U_A)}=-\frac{\lambda_0}{\lambda_1}\frac{\beta^L_B\alpha\tK}{\beta^L_A-\beta^L_B}
                +\frac{-\beta^L_A\alpha\tK}{\beta^L_A-\beta^L_B}+\frac{\alpha\tK\lambda_0}{r+\lambda_0}.
\]
Substituting for (\ref{eq-x*}) these limits and the limits obtained so far, we get the following:
\begin{equation}\label{eq-Sect5-5}
x^*_\infty=\frac{(r-\mu+\lambda_0)\l\{(\lambda_0(\beta^L_A-\zeta^{L,-}_0)\beta^L_B+\lambda_1(\beta^L_B-\zeta^{L,-}_0)\beta^L_A)(r+\lambda_0)
    +\zeta^{L,-}_0(\beta^L_A-\beta^L_B)\lambda_0\lambda_1\r\}\tK}
    {(r+\lambda_0)\l\{(\lambda_0(\beta^L_A-\zeta^{L,-}_0)(\beta^L_B-1)+\lambda_1(\beta^L_B-\zeta^{L,-}_0)(\beta^L_A-1))(r-\mu+\lambda_0)
    +(\zeta^{L,-}_0-1)(\beta^L_A-\beta^L_B)\lambda_0\lambda_1\r\}}
\end{equation}
We can see that $x^*_\infty\geq\displaystyle{\frac{\beta^L_A}{\beta^L_A-1}\tK}\geq\tK$ holds.
In addition, for the case where $\alpha=1$ and $K=0$, we can confirm that the above result coincides with Proposition 1 of \cite{N20}.

%%%%%%%%%%%%%%%%%%%%%%%%%%%%%%%%%%%%%%%%%%%%%%%%%%%%%%%%%%%%%%%%%%%%%%%%%%%%%%%
\subsection{Asymptotic behaviors as $\lambda_0\to\infty$}
As $\lambda_0$ tends to $\infty$, the regime $0$ vanishes, and only the constraint on the random arrival of investment opportunities remains.
In other words, the model converges to the one treated in Dupuis and Wang \cite{DW02}.
In this case, it follows that
\begin{equation}\label{eq-Sect5-2-1}
\l\{\begin{array}{l}
\displaystyle{\lim_{\lambda_0\to\infty}\beta^L_A=\infty, \ \ \ 
\lim_{\lambda_0\to\infty}\beta^U_A=\frac{1}{2}-\frac{\mu}{\sigma^2}-\sqrt{\l(\frac{1}{2}-\frac{\mu}{\sigma^2}\r)^2+\frac{2(\eta+r)}{\sigma^2}},
\ \ \ \lim_{\lambda_0\to\infty}\beta^U_B=-\infty}, \\
\displaystyle{\lim_{\lambda_0\to\infty}a_0,a_1=\frac{\alpha\eta}{r-\mu+\eta}, \ \mbox{ and } \ \lim_{\lambda_0\to\infty}b_0,b_1=-\frac{\alpha\tK\eta}{r+\eta}}. 
\end{array}\r.
\end{equation}
Note that the value of $\beta^L_B$ is independent of $\lambda_0$. We have then $\displaystyle{\lim_{\lambda_0\to\infty}}P^L_A,Q^L_A,P^U_B,Q^U_B=0$, and
\[
\lim_{\lambda_0\to\infty}P^L_B=\alpha, \ \ \ \lim_{\lambda_0\to\infty}Q^L_B=-\alpha\tK, \ \ \
\lim_{\lambda_0\to\infty}P^U_A=\frac{(r-\mu)\alpha}{r-\mu+\eta}, \ \ \ \lim_{\lambda_0\to\infty}Q^U_A=-\frac{r\alpha\tK}{r+\eta}.
\]
Moreover, $G^U_0(\beta^U_A)\sim\eta-\lambda_0$ as $\lambda_0\to\infty$.
By the same way as the previous subsection, we obtain that
\begin{align*}
\lim_{\lambda_0\to\infty}x^*
&= \frac{(r-\mu+\eta)\l\{(\beta^L_B-1)(r+\eta)+(1-\beta^U_A)r+\eta\r\}\tK}{(r+\eta)\l\{(\beta^L_B-1)(r-\mu+\eta)+(1-\beta^U_A)(r-\mu)\r\}} \\
&= \frac{(r-\mu+\eta)((r+\eta)\beta^L_B-r\beta^U_A)\tK}{(r+\eta)((r-\mu+\eta)\beta^L_B-(r-\mu)\beta^U_A-\eta)}(=:x^*_\infty),
\end{align*}
\[
\lim_{\lambda_0\to\infty}v_1(x)= \l\{\begin{array}{l}
          \displaystyle{\alpha(x^*_\infty-\tK)\l(\frac{x}{x^*_\infty}\r)^{\beta^L_B}}, \ \ \ 0<x<x^*_\infty, \\
          \displaystyle{\l(\frac{(r-\mu)\alpha x^*_\infty}{r-\mu+\eta}-\frac{r\alpha\tK}{r+\eta}\r)\l(\frac{x}{x^*_\infty}\r)^{\beta^U_A}
          +\frac{\alpha\eta}{r-\mu+\eta}x-\frac{\alpha\tK\eta}{r+\eta}}, \ \ \ x>x^*_\infty,
        \end{array}\r.
\]
and $\displaystyle{\lim_{\lambda_0\to\infty}v_0(x)=\lim_{\lambda_0\to\infty}v_1(x)}$ for any $x>0$, where $\beta^U_A$ is the limit given in (\ref{eq-Sect5-2-1}).
As seen in \cite{DW02}, we can prove that $\displaystyle{x^*_\infty\geq \frac{r(r-\mu+\eta)}{(r-\mu)(r+\eta)}\tK}\geq\tK$ holds and
the boundary conditions (\ref{VI-*L}) and (\ref{VI-*U}) are satisfied.
When $\alpha=0$ and $I=0$, the result in this subsection is consistent with \cite{DW02}.

%%%%%%%%%%%%%%%%%%%%%%%%%%%%%%%%%%%%%%%%%%%%%%%%%%%%%%%%%%%%%%%%%%%%%%%%%%%%%%%
%
% Section 6
%
%%%%%%%%%%%%%%%%%%%%%%%%%%%%%%%%%%%%%%%%%%%%%%%%%%%%%%%%%%%%%%%%%%%%%%%%%%%%%%%
\setcounter{equation}{0}
\section{Conclusions}
We considered a two-state regime-switching model and discussed the optimal stopping problem defined by (\ref{def-vi}) under two constraints on stopping:
the random arrival of investment opportunities and the regime constraint. 
Under the assumption that the boundary conditions (\ref{VI-*L}) and (\ref{VI-*U}) are satisfied, we showed that an optimal stopping time exists as a threshold type.
In addition, we derived expressions of the value functions $v_i, i=0,1$, and the optimal threshold $x^*$, which include solutions to quartic equations, but can be easily computed numerically.
Asymptotic behaviors of $v_i, i=0,1$, and $x^*$ are also discussed.
On the other hand, the assumption of the boundary conditions might be redundant, as mentioned in Remark \ref{rem3-1}.
Thus, it is significant as future work to show that the boundary conditions (\ref{VI-*L}) and (\ref{VI-*U}) are always satisfied using, e.g.,
a PDE approach discussed in Bensoussan et al. \cite{BYY12}.

%%%%%%%%%%%%%%%%%%%%%%%%%%%%%%%%%%%%%%%%%%%%%%%%%%%%%%%%%%%%%%%%%%%%%%%%%%%%%%%
\begin{center}
{\bf Acknowledgments}
\end{center}
Takuji Arai gratefully acknowledges the financial support of the MEXT Grant in Aid for Scientific Research (C) No.18K03422.

%%%%%%%%%%%%%%%%%%%%%%%%%%%%%%%%%%%%%%%%%%%%%%%%%%%%%%%%%%%%%%%%%%%%%%%%%%%%%%%

\end{document}